\documentclass{gtpart}     
\usepackage{graphicx}  
\usepackage[all]{xy}
\usepackage{amscd}
\usepackage{psfrag}
\usepackage{amsmath}
\usepackage{amssymb}
\usepackage{epstopdf}

\title{Epimorphisms and Boundary Slopes of 2-Bridge Knots}

\author{Jim Hoste}
\givenname{Jim}
\surname{Hoste}
\address{Pitzer College}
\email{jhoste@pitzer.edu}
\urladdr{}

\author{Patrick D. Shanahan}
\givenname{Patrick}
\surname{Shanahan}
\address{Loyola Marymount University}
\email{pshanahan@lmu.edu}
\urladdr{}


\arxivreference{}
\arxivpassword{}

\volumenumber{}
\issuenumber{}
\publicationyear{}
\papernumber{}
\startpage{}
\endpage{}
\doi{}
\MR{}
\Zbl{}
\received{}
\revised{}
\accepted{}
\published{}
\publishedonline{}
\proposed{}
\seconded{}
\corresponding{}
\editor{}
\version{}

\newtheorem{theorem}{Theorem}[section]
\newtheorem{lemma}[theorem]{Lemma}
\newtheorem{corollary}[theorem]{Corollary}
\newtheorem{question}[theorem]{Question}
\newtheorem{conjecture}[theorem]{Conjecture}

\begin{document}

\begin{abstract}    
In this article we study a partial ordering on 
knots in $S^3$ where $K_1 \ge K_2$ if there is an epimorphism 
from the knot group of $K_1$ onto the knot group of $K_2$ which 
preserves peripheral structure.  If $K_1$ is a 2-bridge knot and $K_1 
\ge K_2$, then it is known  that $K_2$ must also be 2-bridge. Furthermore, Ohtsuki, Riley, and Sakuma give a construction which,  for a given 2-bridge knot $K_{p/q}$, produces infinitely  2-bridge knots $K_{p'/q'}$ with  $K_{p'/q'}  \ge K_{p/q}$.  After characterizing all 2-bridge knots with 4 or less distinct boundary slopes, we use this to prove that in any such pair,  $K_{p'/q'}$ is either a torus knot or has 5 or more distinct boundary slopes.  We also prove that 2-bridge knots with exactly 3 distinct boundary slopes are minimal with respect to the partial ordering. This result provides some evidence for the conjecture that all pairs of 2-bridge knots with $K_{p'/q'} \ge K_{p/q}$ arise from the Ohtsuki-Riley-Sakuma construction. 
\end{abstract}

\maketitle


\section{Introduction}

Many interesting problems in knot theory involve the question of when one knot complement can be mapped onto another, or at the algebraic level, when the fundamental group of one knot complement can be mapped onto the fundamental group of another.  For example, Simon's Conjecture asserts that the fundamental group of any knot in $S^3$ can surject onto only finitely many knot groups.  This conjecture has been established for 2-bridge knots by Boileau, Boyer, Reid and Wang \cite{BBRW:2009}.  On the other hand, given any 2-bridge knot $K_{p/q}$, Ohtsuki, Riley, and Sakuma \cite{ORS} show how to construct infinitely many 2-bridge knots $K_{p'/q'}$ whose complements map onto the complement of $K_{p/q}$. Their construction induces an epimorphism $\phi: \pi_1(S^3-K_{p'/q'}) \to \pi_1(S^3-K_{p/q})$ of knot groups which {\it preserves peripheral structure}, that is, $\phi$ takes the subgroup generated by a longitude and meridian of $K_{p'/q'}$ into a conjugate of the subgroup generated by a longitude and 
meridian of $K_{p/q}$.

For any pair of knots, not just 2-bridge knots,  Silver and Whitten define a partial ordering $\ge$ on 
knots in $S^3$ by declaring that $K_1\ge K_2$ if there exists 
an epimorhism $\phi: \pi_1(S^3 - K_1) \to \pi_1(S^3 - 
K_2)$ which preserves peripheral structure  \cite{SW:2006}. Thus, for a given 2-bridge knot $K_{p/q}$, the Ohtsuki-Riley-Sakuma (ORS) construction gives infinitely many 2-bridge knots $K_{p'/q'}$ such that $K_{p'/q'} \ge K_{p/q}$. This leads naturally to the following question posed in \cite{ORS}.

\begin{question}\label{ORS question}(Ohtsuki-Riley-Sakuma) Is every pair of 2-bridge knots  $(K_{p'/q'}, K_{p/q})$ with $K_{p'/q'} \ge K_{p/q}$ given by the ORS construction?
\end{question}

Gonz\'{a}lez-Acu\~{n}a and Ram\'{i}rez \cite{GR:2001} give an affirmative answer to Question~\ref{ORS question} 
in the case where $K_{p/q}$ is a 2-bridge torus knot.  In this paper we give additional evidence that the answer to Question~\ref{ORS question} is yes in general. Numerous computations with non-torus knot ORS pairs $(K_{p'/q'}, K_{p/q})$ with $K_{p'/q'} \ge K_{p/q}$ suggest that in any such pair the larger knot $K_{p'/q'}$ must be sufficiently complex with respect to a variety of knot invariants such as: crossing number, degree of the Alexander polynomial, degree of the A-polynomial, number of distinct boundary slopes, etc. Thus, if the answer to Question~\ref{ORS question} is yes, then knots with ``small'' complexity should be minimal.  Here, a knot $K_1$ is said to be {\em minimal} with respect to the Silver--Whitten partial ordering if $K_1\ge K_2$ implies that either $K_2=K_1$, $K_2$ is the mirror image of $K_1$, or $K_2$ is the unknot. Throughout the remainder of this paper, we will consider two knots $K_1$ and $K_2$ to be equivalent, $K_1 \equiv K_2$,  if either $K_2=K_1$ or $K_2$ is the mirror image of $K_1$.

One of the themes of this paper will be to look at the set of distinct boundary slopes of a knot. These are easily computed for 2-bridge knots due to the classification of Hatcher and Thurston \cite{HT:1985}. 
Our first main result gives a lower bound on the number of distinct boundary slopes for the larger knot in an ORS pair.

\setcounter{section}{4}
\begin{theorem}
If $K_{p'/q'} \ge K_{p/q}$ is a nontrivial ORS pair, then either
\begin{enumerate}
 \item[{\it i}.] $K_{p'/q'} $ and $K_{p/q}$ are both torus knots and $K_{p'/q'} $ has precisely two distinct boundary slopes, or
 \item [{\it ii}.]$K_{p'/q'} $ has at least five distinct boundary slopes.
 \end{enumerate}
 \label{thm1}
 \end{theorem}
\setcounter{section}{1}

Our second main result establishes the minimality of 2-bridge knots with exactly three boundary slopes.

\setcounter{section}{5}
\setcounter{theorem}{0}
\begin{theorem} 
If $K_{p/q}$ is a 2-bridge knot with exactly 3 distinct boundary slopes, then $K_{p/q}$ is minimal with respect to the Silver--Whitten partial order.
\label{thm2}
\end{theorem}
\setcounter{section}{1}

Notice that if the answer to Question~\ref{ORS question} is yes, then Theorem~\ref{thm1} implies Theorem~\ref{thm2}. Thus, we may think of Theorem~\ref{thm2} as evidence in support of an affirmative answer to Question~\ref{ORS question}. It is also worth noting that Theorem~\ref{thm1} is sharp. That is, there exist non-minimal 2-bridge knots with 5 distinct slopes. We provide an example of this in Section~4.

This project began as joint work of the first author and two undergraduate students, T\'omasz Przytycki and Rebecca Nachison, in the summer of 2007 at the Claremont Colleges REU program. Thanks are due to T\'omasz and Rebecca for working out most of Theorem~\ref{list of knots with given number of slopes}. We thank the National Science Foundation, the Claremont Colleges, and Pitzer College in particular, for their support of the REU program.  The second author would also like to thank the Claremont Colleges for their hospitality and Pitzer College for its support during the completion of this work.   Finally, thanks are due to Alan Reid for helpful conversations regarding character varieties.

\section{Preliminaries}    
We begin with some notation for 2-bridge knots. Recall that each 2-bridge knot  corresponds to a relatively prime pair of integers $p$ and $q$ with $q$ odd. We denote the knot as $K_{p/q}$. Furthermore, $K_{p/q}$ and $K_{p'/q'}$ are ambient isotopic as unoriented knots if and only if $q'=q$ and $p' \equiv p^{\pm 1}\  (\mbox{mod } q)$ (see \cite{BZ:2003} for details).  As mentioned earlier, we will not distinguish between a knot $K_{p/q}$ and its mirror image $K_{-p/q}$.  Therefore, in this paper we consider two 2-bridge knots $K_{p/q}$ and $K_{p'/q'}$ to be equivalent if and only if $q'=q$ and either $p' \equiv p^{\pm 1}\  (\mbox{mod } q)$ or $p' \equiv -p^{\pm 1}\  (\mbox{mod } q)$. Because of this classification, it is sometimes convenient to assume that $0<p<q$. 

The theory of 2-bridge knots is closely tied to continued fractions. We adopt the convention used in \cite{ORS} and define $p/q=r+[b_1, b_2, \dots, b_n]$ as the continued fraction
 
$$p/q=r+[b_1, b_2, \dots, b_n]=r+\frac{1}{b_1+\displaystyle \frac{1}
{
\begin{array}{ccc}
b_2+&&\\
&\ddots&\\
&&+\displaystyle \frac{1}{b_n}\\
 \end{array}
 }
 }$$

It is well known that if $0<p<q$ then we may write $p/q$ uniquely as $p/q=0+[b_1, b_2, \dots , b_n]$ where each $b_i>0$ and $b_n>1$. Furthermore, we may also assume that $b_1>1$ for the following reason. If $b_1=1$, reversing the order of the partial quotients will give $p'/q=[b_n, \dots,b_3, b_2, 1]=[b_n, \dots, b_3, b_2+1]$ with $0<p'<q$ and $p p'\equiv (-1)^{n+1}\ \text{ (mod $q$)}$. Thus, for any 2-bridge knot $K$, there exists a continued fraction $0+[b_1, b_2, \dots, b_n]$ with $b_i>0, \ b_1>1$, and $b_n>1$ representing $K$ (or its mirror image). In what follows, we will refer to such a continued fraction as {\it strongly positive}. Moreover, we will also refer to a vector $(b_1, b_2, \dots, b_n)$ of positive integers with $b_1 > 1$ and $b_n >1$ as strongly positive. Finally, the negation of such a vector we call {\it strongly negative}.

\begin{figure}[t]
\psfrag{a}{\Huge$b_1$}
\psfrag{b}{\Huge$b_2$}
\psfrag{c}{\Huge$b_3$}
\psfrag{d}{\Huge$b_n$}
     \begin{center}
     \leavevmode
     \scalebox{.4}{\includegraphics{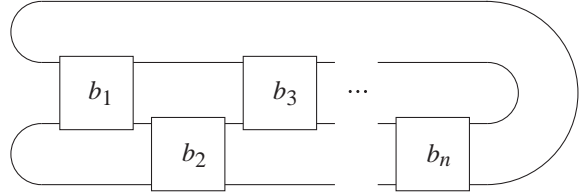}}
     \end{center}
\caption{A 4-plat diagram for $K_{p/q}$ where $p/q = r+[b_1, b_2, \dots, b_n]$ ($n$ even).}
\label{4plat fig}
\end{figure}

If $p/q=r+[b_1, b_2, \dots, b_n]$, then $K_{p/q}$ can be represented by the standard 4-plat diagram shown in 
Figure~\ref{4plat fig}. In each box there are $b_i$ crossings.  If $i$ is odd, then the crossings are right-handed if $b_i >0$ and left-handed otherwise.  For $i$ even the convention reverses, with $b_i > 0$ corresponding to left-handed crossings and $b_i < 0$ for right-handed.  Thus, the 4-plat diagram is alternating precisely when all $b_i$ have the same sign. Notice also that there are two possible ways to ``close up'' the strands at the right end of the braid in Figure~\ref{4plat fig} depending on whether $n$ is even or odd. Figure~\ref{4plat fig} depicts a 4-plat with $n$ even. If $n$ is odd, strings 1 and 2, and strings 3 and 4, connect pairwise at the right end of the braid just as they do at the left end. With this notation, the fraction $1/3=[3]$ gives the right-hand trefoil knot while $2/5=[2,2]$ gives the figure eight knot.

Now suppose that $\phi$ is a  homomorphism between knot groups that preserves peripheral structure. Since all nontrivial knots in $S^3$ have Property P \cite{KM:2004}, there are certain restrictions placed on the images of a meridian and longitude of $K_1$.

\begin{theorem} Let $K_1$ and $K_2$ be nontrivial knots with $K_1\ge 
K_2$. Then there exists an epimorphism $\phi: \pi_1(
S^3-K_1) \to \pi_1(S^3-K_2)$ and an integer $d$ such that
\begin{eqnarray}
\phi(\mu_1)&=&\mu_2, \mbox{ and}\\
\phi(\lambda_1)&=&\lambda_2^d
\end{eqnarray}
for some choice $\{ \mu_i, \lambda_i \}$ of meridian-(preferred)\,longitude pairs 
for $K_i$.
\label{abcdcriteria}
\end{theorem}
\noindent {\em Proof.}  Since $\phi$ preserves peripheral structure we 
may assume (after conjugation, if necessary) that
\begin{eqnarray*}
\phi(\mu_1)&=&\mu_2^a \lambda_2^b, \mbox{ and}\\
\phi(\lambda_1)&=&\mu_2^c \lambda_2^d
\end{eqnarray*}
for some integers $a, b, c$ and $d$. If we consider the 
abelianization map $ab:\pi_1(S^3-K_2) \rightarrow \mathbb Z$, 
then $ab(\phi (\lambda_1)) = 0$ since $\lambda_1$ is in the 
commutator subgroup.  However,   $ab(\phi (\lambda_1)) = ab ( \mu_2^c 
\lambda_2^d) = c$. This shows that $c=0$. Moreover, $ab(\phi (\mu_1)) 
= ab(\mu_2^a \lambda_2^b)=a$ must generate $\mathbb Z$ because $ab 
\circ \phi$ is onto and $\pi_1(S^3-K_1)$ is generated by 
conjugates of $\mu_1$.  Therefore, $a = \pm 1$. By replacing $\mu_2$ 
with its inverse if necessary, we may assume that $a=1$.

Now consider $1/b$ Dehn filling on $K_2$. This filling kills the 
image of $\mu_1$ and any of its conjugates, and so kills the entire 
group, since $\phi$ is onto. Thus the manifold we obtain after filling is simply connected. 
Since $K_2$ is a nontrivial knot with Property P, we must have done 
trivial surgery. Therefore, $b=0$.
\hfill $\Box$

Because of Theorem~\ref{abcdcriteria}, it makes sense to write $$K_1\ge_dK_2$$ whenever $K_1\ge K_2$ and there exists an epimorphism $\phi$ as in Theorem~\ref{abcdcriteria} with $\phi(\lambda_1)=\lambda_2^d$.  In general, $d$ can take on any integer value. In particular, given the knot $K_2$ and any integer $d$, there exists a knot $K_1$ and epimorphism $\phi$ realizing $d$. The reader is referred to \cite{JL:1989} for entr\'e into this subject.  On the other hand, if $K_1$ and  $K_2$ are 2-bridge knots, then  $d$ must be odd.

\begin{theorem}
\label{2bridge means d is odd} 
If $K_{p'/q'}\ge_d K_{p/q}$  where $K_{p'/q'}$ and $K_{p/q}$ are 2-bridge knots and $K_{p/q}$ is nontrivial, then $d$ is odd. In particular, the image $\phi(\lambda_1)$ under the epimorphism $\phi$ cannot be trivial.
\end{theorem}

\noindent {\em Proof.} Suppose $\phi:\pi_1(S^3-K_{p'/q'}) \to 
\pi_1(S^3-K_{p/q})$  is an epimorphism taking  $\lambda_1$ to $\lambda_2^d$ and let $\rho$ be an irreducible  
parabolic representation of $\pi_1(S^3-K_{p/q})$ into $SL(2, \mathbb C)$ (such a representation exists by Riley \cite{R:1984}). 
Composing $\rho$ with $\phi$ gives such a representation for $K_{p'/q'}$. 
It follows from \cite{HS:2004} that any irreducible parabolic representation of a 
2-bridge knot into $SL(2,\mathbb C)$  must take the longitude into the conjugacy class of an 
upper triangular matrix with diagonal entries of $-1$. Hence $\phi$ 
cannot take $\lambda_1$ to an even power of $\lambda_2$.
\hfill $\Box$

The following theorem summarizes some additional properties of the partial ordering which we will refer to in this article.  The first two parts are contained in \cite{SW:2006}, the third comes from \cite{BBRW:2009}, and the fourth appears in \cite{BB:2007}.

\begin{theorem}
\label{properties of partial order} 
(Silver-Whitten, Boileau-Boyer, Boileau-Boyer-Reid-Wang) Suppose that $K_1 \ge_d K_2$.
\begin{enumerate}
\item[{\it i}.] If $\Delta^{(i)}_K(t)$  is the $i$-th Alexander polynomial of $K$, then $\Delta^{(i)}_{ K_2}(t)$ divides $\Delta^{(i)}_{K_1}(t)$.
\item[{\it ii}.] If $A_K(L,M)$ is the A-polynomial of $K$, then 
$A_{K_2}(L,M)$ divides $(L^d-1)A_{K_1}(L^d,M)$.
\item[{\it iii}.] If $X(K)$ is the $SL(2, \mathbb C)$ character variety of $K$, then the induced map $\phi^* : X(K_2) \rightarrow X(K_1)$ is an injective, algebraic  and closed (in the Zariski topology) mapping.
\item[{\it iv}.] If $K_1$ is a 2-bridge knot, then $K_2$ is a 2-bridge knot. Futhermore, if $K_1=K_{p'/q'}$ and $K_2=K_{p/q}$, then either $K_1 \equiv K_2$ or  $q'=kq$ with $k>1$.
\end{enumerate}
\end{theorem}

Notice that the fact that $q$ divides $q'$ in Theorem~\ref{properties of partial order}({\em iv}) is easily derived from ({\em i}) and the fact that $q=| \Delta_{K_{p/q}}(-1) |$ is the determinant of $K_{p/q}$.  

Because of the relationship between the A-polynomial and boundary slopes of a knot,  Theorem~\ref{properties of partial order}({\it ii}) enables us to use boundary slopes as a tool to study whether one knot is greater than another.  Recall that a {\it boundary slope} $r=\frac{a}{b}$ of a knot $K$ is an element of $\mathbb Q \cup \infty$ such that $K$ contains a properly embedded,   incompressible, boundary incompressible surface $S$ in its exterior whose boundary is the curve $ a \mu + b \lambda$ (or multiple copies of this curve). Here $\{\mu, \lambda\}$ is a meridian-(preferred) longitude pair for $K$.  The {\it Newton polygon} of the A-polynomial $A_K(L,M)=\sum a_{i\,j}L^i M^j$ is the convex hull of the set of points $\{(i,j) \in \mathbb R^2\, |\, a_{i\, j} \ne 0\}$.  If $r$ is the slope of a side of the Newton polygon of $A_K(L,M)$, then it is proven in \cite{CCGLS:1994}, that $r$ is a boundary slope of $K$. A boundary slope of $K$ which appears as a slope of the Newton polygon is called {\it strongly detected}.  We shall tacitly assume that the A-polynomial always includes the component $L-1$ corresponding to abelian representations so that $0$ is always a strongly detected slope.

\begin{corollary} Suppose $K_1\ge_d K_2$, $d\ne 0$,  and 
that $r$ is a strongly detected boundary slope of $K_2$. Then $dr$ is 
a strongly detected boundary slope of   $K_1$. 
\end{corollary}
\noindent {\em Proof.} We establish this result using the following well-known property of Newton polygons.  For any polynomial $P(x,y)$ let $S_P$ denote the set of slopes of the sides of the Newton polygon of $P$.   Then for all polynomials $P(x,y)$ and $Q(x,y)$, $S_P \cup S_Q = S_{PQ}$ (a proof may be found in \cite{HS:2004}).  Now suppose that $r$ is a strongly detected boundary slope of $K_2$.  Then $r$ is the slope of a side of the Newton polygon of $A_{K_2}(L,M)$. Since $A_{K_2}(L,M)$ divides 
$(L^d-1)A_{K_1}(L^d,M)$ it follows that $r=0$ or $r$ is a slope of a side of
the Newton polygon of $A_{K_1}(L^d,M)$. If $r=0$, then $dr=0$ is a strongly detected slope of $K_1$. Now  assume that 
$r$ is a slope of a side of the Newton polygon of $A_{K_1}(L^d,M)$.  The Newton polygon of $A_{K_1}(L^d,M)$ is obtained from the Newton polygon of $A_{K_1}(L,M)$ by simply replacing every vertex $(i,j)$ with $(i/d,j)$.
Consequently, $dr$
is the slope of a side of the Newton polygon of $A_{K_1}(L,M)$, and therefore, is a strongly detected boundary slope of $K_1$. \hfill $\Box$

By work of Ohtsuki \cite{O:1994}, all boundary slopes of 2-bridge knots are strongly detected. Thus, for 2-bridge knots we have a stronger result.

\begin{corollary} 
\label{boundary slopes are subsets}
Suppose $K_1\ge_d K_2$ and that 
$K_1$ is a 2-bridge knot. If $\{r_1, r_2, 
\dots, r_m\}$ is the set of boundary slopes of $K_2$, then  $\{dr_1, 
dr_2, \dots, dr_m\}$ is a subset of the set of boundary slopes of 
$K_1$.
\end{corollary}

For any knot $K$, the {\it diameter} of its set of boundary slopes is the maximum difference between any two slopes. We denote the diameter by $diam(K)$. For 2-bridge knots, in fact for all alternating Montesinos knots, $cr(K)=diam(K)/2$, where $cr(K)$ is the crossing number of $K$. See \cite{HS:2007}, \cite{IM:2008} and \cite{MMR:2008}. This leads to the following results.

\begin{corollary} 
\label{crossing number}
Suppose $K_1$ is a 2-bridge knot and $K_1\ge_d K_2$. 
Then $cr(K_1)\ge |d|cr(K_2)\ge cr(K_2)$.
\end{corollary}

Using Theorem~\ref{properties of partial order}({\it iii}) and Corollary~\ref{crossing number} it is now easy to show the following. (This result was first proven by Boileau, Rubinstein, and Wang \cite{BRW:2005}.)
 
\begin{corollary} A 2-bridge knot can only be greater than or equal 
to finitely many other knots.
\end{corollary}

\section{2-bridge knots with four or less boundary slopes}

In preparation for proving Theorem~\ref{thm1}, we determine in this section those 2-bridge knots with small numbers of distinct boundary slopes. In \cite{MMR:2008}, the authors determine necessary but not sufficient conditions on $p$ and $q$ so that $K_{p/q}$ has 4 or less distinct boundary slopes.  In order to exploit Corollary~\ref{boundary slopes are subsets} we require a stronger result than what appears in \cite{MMR:2008}.  In particular, we need both a complete classification of 2-bridge knots with 4 or less boundary slopes as well as an explicit description of the associated slope sets.

We assume the reader is familiar with Hatcher and Thurston's paper \cite{HT:1985} where a method for computing the boundary slopes of $K_{p/q}$ is given. We will also make use  of a method equivalent to theirs which is described in \cite{HS:2007}. Rather than describing the methodology in complete detail, we provide a brief description illustrated by an example.

Recall that the Farey graph may be thought of as the edges of the ideal modular tessellation of $\mathbb H^2$ obtained by starting with the ideal triangle whose vertices are $1/0,\ 0/1$, and $1/1$ and then reflecting it in all possible ways across its edges. The ideal vertices of this tessellation are $\mathbb Q \cup \{\infty\}$. Two vertices  $a/b$ and $c/d$ are joined by an edge precisely when $\det \begin{pmatrix}a&c\\b&d\end{pmatrix}=\pm 1$ and we call this value the {\it determinant} of the directed edge from $a/b$ to $c/d$. If $a/b$ and $c/d$ are joined by an edge, then so are $a/b$ and $(a+c)/(b+c)$ as well as $c/d$ and $(a+c)/(b+d)$. In this case the two fractions  $a/b$ and $c/d$, together with their {\it mediant} $(a+c)/(b+c)$, form the vertices of an ideal triangle in the tessellation. Finally, given any continued fraction expansion of $p/q$,  the associated sequence of convergents defines a path in the Farey graph from $1/0$ to $p/q$. The number of triangles through which the path ``turns'' at each vertex is equal to the corresponding partial quotient.  

\begin{figure}
\psfrag{a}{\huge $\frac{1}{0}$}
\psfrag{b}{\huge $\frac{1}{1}$}
\psfrag{c}{\huge $\frac{1}{2}$}
\psfrag{d}{\huge $\frac{3}{7}$}
\psfrag{e}{\huge $\frac{5}{12}$}
\psfrag{f}{\huge $\frac{7}{17}$}
\psfrag{g}{\huge $\frac{0}{1}$}
\psfrag{h}{\huge $\frac{1}{3}$}
\psfrag{i}{\huge $\frac{2}{5}$}
    \begin{center}
    \leavevmode
    \scalebox{.5}{\includegraphics{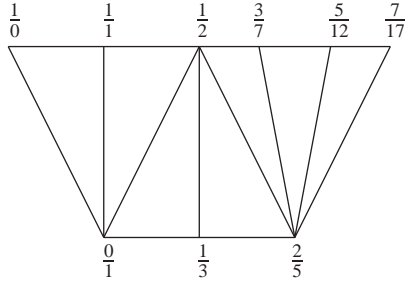}}
    \end{center}
\caption{The portion of the Farey graph between $1/0$ and $7/17$.}
\label{7/17}
\end{figure}

For example, consider $7/17=0+[2,2,3]$. We picture the relevant portion of the Farey graph in Figure~\ref{7/17}. The sequence of convergents is $\{0, 1/2, 2/5, 7/17\}$. The partial quotients are $2,2,3$ and the path turns through 2 triangles on the left at $0$, 2 triangles on the right at $1/2$ and 3 triangles on the left at $2/5$. If all the partial quotients are positive, as in this example, then the turning at each vertex alternates between triangles on the left and right. Negative partial quotients correspond to turning in the opposite direction. Again, the reader is urged to consult \cite{HT:1985} and \cite{HS:2007}.

For any $p/q$, there are infinitely many edge paths in the Farey graph that begin at $1/0$ and end at $p/q$, but only finitely many that are {\it minimal}. A path is minimal if it never backtracks and if it never traverses two edges of the same triangle. The minimal paths for $K_{7/17}$ and their corresponding continued fractions are shown in the first two columns of Table~\ref{slope data 7/17}.  In general, only one of the minimal paths is {\it even}, that is, all the partial quotients are even. In Table~\ref{slope data 7/17}, the first path is even. According to Hatcher and Thurston, each minimal path determines a boundary slope. The slope can be computed from the path $\gamma$ by the formula
$$-2(m(\gamma)-m(\gamma_{even}))$$
where $\gamma_{even}$ is the even path and $m(\gamma)$ is the sum of the determinants of each edge of the path, excluding the first edge. Thus, the even path always gives a slope of zero.  For any path $\gamma$, we call $m(\gamma)$ the {\it unadjusted} slope. The last two columns of Table~\ref{slope data 7/17} give the unadjusted slopes and the boundary slopes of $K_{7/17}$. 

\begin{table}
\begin{center}
\caption{Boundary slope data for $K_{7/17}$}
\vskip .1 in
\begin{tabular}{l|l|l|l}
minimal path $\gamma$ & fraction & $m(\gamma)$ & slope\\
\hline
 $\{1/0, 1/1, 1/2, 3/7, 5/12, 7/17\}$&$1+[-2,4,-2,2]$& 4 & 0\\
 $ \{1/0, 0/1, 1/2, 3/7, 5/12, 7/17\}$&$0+[2,3,-2]$& 2 & 4\\
  $ \{1/0, 1/1, 1/2, 2/5, 7/17\}$&$1+[-2,3,3]$& 1 & 6\\
   $\{1/0, 0/1, 1/2, 2/5, 7/17\}$&$ 0+[2,2,3]$& -1 & 10\\
$\{1/0, 1/3, 2/5, 7/17\}$&$0+[3,-2,4]$& -3 & 14
\end{tabular}
\label{slope data 7/17}
\end{center}
\end{table}

For any $p/q$, if $\gamma_1$ and $\gamma_2$ are two paths from $1/0$ to $p/q$, then we can move from one path to the other by a sequence of {\it triangle moves}. By a triangle move, we mean replacing two consecutive edges that lie in a single triangle with the third edge of the triangle, or vice-versa. If $\gamma_1$ is changed to $\gamma_2$ by a single triangle move, and furthermore, the triangle lies on the right of $\gamma_1$, then we call the move a {\it right} triangle move. Right triangle moves define a partial order on the set of paths: we say that $\gamma_1>\gamma_2$ if $\gamma_1$ can be changed to $\gamma_2$ by a sequence of right triangle moves. Furthermore, it is shown in \cite{HS:2007} that right triangle moves strictly decrease the unadjusted slope of the path, thus, $\gamma_1 > \gamma_2$ implies that $m(\gamma_1) > m(\gamma_2)$.  The {\it upper} and {\it lower} minimal paths can now be defined as ones that are either maximal or minimal with respect to this partial order, respectively. These paths are well-defined because there is a unique minimal path with no triangles on its left and a unique minimal path with no triangles on its right.  In Table~\ref{slope data 7/17}, the first path is the upper path and the last path is the lower path.   

The following lemma gives a lower bound on the number of distinct boundary slopes.

\begin{lemma}
Suppose $p/q=r+[a_1, a_2, \dots, a_m]$ is a strongly positive continued fraction. Then the number of distinct  boundary slopes of $K_{p/q}$ is at least $2+\lfloor m/2\rfloor$ and at most $f_{m+1}$ where $f_m$ denotes the $m$-th Fibonacci number (assuming $f_0=0$ and $f_1=1$).
\label{bounds on number of paths}
\end{lemma}

\noindent {\em Proof.}  For the lower bound we will create a sequence of minimal paths $\gamma_\text{upper}$, $\gamma_1$, $\dots$, $\gamma_{\lfloor m/2 \rfloor}$, $\gamma_\text{lower}$ with strictly decreasing unadjusted slopes. Let $\gamma_\text{upper}$ be the upper path. Next, let $\gamma_1$ be the path that goes from $1/0$ to $r/1$ to $r+[a_1]$ and then continues with the rest of the upper path. The path $\gamma_1$ is obtained from $\gamma_\text{upper}$ by a sequence of $a_1$ right triangle moves.  Now let $\gamma_2$ begin  $\{1/0, r/1, r+[a_1, a_2], r+[a_1, a_2, a_2]\}$ and then continue with the rest of the upper path. Notice that $\gamma_2$ is obtained from $\gamma_1$ by a sequence of $a_2 + a_3$ right triangle moves. We continue in this manner to generate each path $\gamma_i$.  The fact that both $a_1>1$ and $a_m>1$ guarantee that all of these paths are minimal. Finally, their unadjusted slopes are strictly decreasing since each path is obtained from the previous one by performing right triangle moves. This establishes the lower bound. It is shown in \cite{HT:1985} that $f_m$ is an upper bound on the number of minimal paths, each of which might give a distinct boundary slope.
\hfill $\Box$

As the number $m$ of partial quotients grows, it becomes tedious to characterize all knots with $m$ partial quotients and exactly $k$ distinct boundary slopes. Therefore, in the following theorem, we only classify 2-bridge knots with up to four distinct boundary slopes. Furthermore, we provide the set of slopes for only one knot from each chiral pair. (Remember that switching from a knot to its mirror image will negate the set of boundary slopes.)

\begin{theorem} 
\label{list of knots with given number of slopes} Let $K_{p/q}$ be a 2-bridge knot with 2, 3, or 4, distinct boundary slopes. Then a  strongly positive continued fraction representing $K_{p/q}$, or its mirror image, together with the associated slope set, is given in the following table.

\begin{center}
\begin{tabular}{|l|l|l|}
\hline
& continued fraction & slope set \\

& &*=multiplicity 2, **=multiplicity 3\\
\hline

$2$ &  $[a_1], a_1$ odd& $\{0, 2a_1\}$\\
\hline

$3_i$ &  $[a_1, a_2]$,  $a_1$ even, $a_2$ even& $\{-2a_1, 0, 2 a_2\}$\\

$3_{ii}$ &  $[a_1, a_2]$, $a_1$ even, $a_2$ odd&  $\{0, 2a_1, 2a_1+2a_2\}$\\

$3_{iii}$&  $[a_1, 1, a_1], a_1$ odd& $\{-4a_1-2,-2a_1-2^*,0\}$\\
\hline

$4_i$ &  $[a_1, a_2, a_1]$, $a_1$ odd, $a_2$ odd, $a_2 > 1$ &$\{-4a_1-2a_2, -2a_1-2a_2^*, -2a_2, 0\}$\\

$4_{ii}$ &  $[a_1, 1, a_3]$,  $a_1\ne a_3$, $a_1$ odd, $a_3$ odd  & $\{-2a_1-2a_3-2,-2a_3-2,-2a_1-2,0\}$\\

$4_{iii}$ &  $[a_1, 1, a_3]$, $a_1 \ne a_3$, $a_1$ even, $a_3$ odd  &$\{-2a_1,0,-2a_1+2a_3,2a_3+2\}$\\

$4_{iv}$&  $[a_1,1,a_1,a_1+1],  a_1$ even &$\{-2a_1, 0^*, 2a_1+2^*, 4a_1+4\}$\\

$4_{v}$&  $[2,1,1,1,2]=8/21$&$\{-8, -4^*, 0^{**}, 6\}$\\ \hline
\end{tabular}
\end{center}
\end{theorem}

\begin{figure}[t]
\psfrag{a}{\Huge $a_1$}
\psfrag{b}{\Huge $a_2$}
\psfrag{c}{\Huge $a_3$}
\psfrag{d}{\Large $1/0$}
\psfrag{e}{\Large $0/1$}
\psfrag{f}{\Large $[a_1]$}
\psfrag{g}{\Large $[a_1,a_2]$}
\psfrag{h}{\Large $[a_1,a_2,a_3]$}
\psfrag{i}{\Large $\gamma_\text{upper}$}
\psfrag{j}{\Large $m=a_1+a_3-1$}
\psfrag{k}{\Large $\gamma_1$}
\psfrag{l}{\Large $m=a_3-1$}
\psfrag{m}{\Large $\gamma_2$}
\psfrag{n}{\Large $m=a_1-1$}
\psfrag{o}{\Large $\gamma_3$}
\psfrag{p}{\Large $m=1$}
\psfrag{q}{\Large $\gamma_\text{lower}$}
\psfrag{r}{\Large $m=-a_2-1$}
     \begin{center}
     \leavevmode
     \scalebox{.4}{\includegraphics{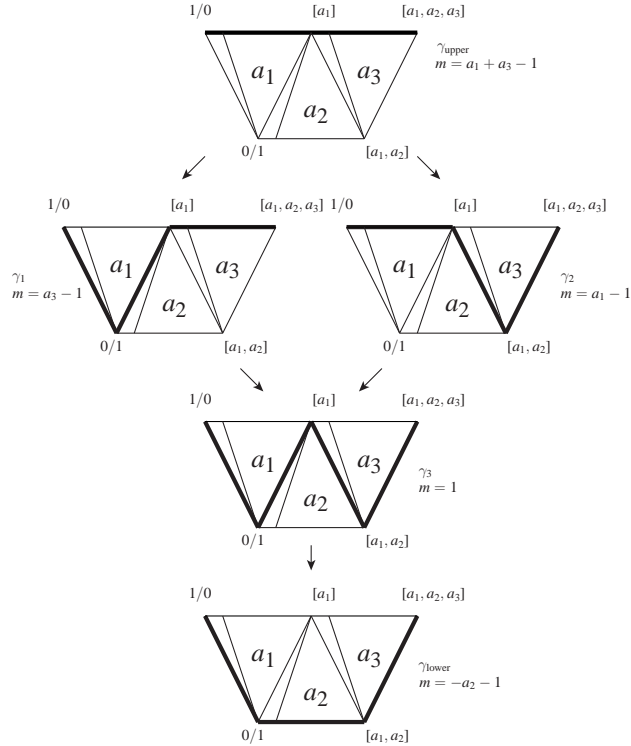}}
     \end{center}
\caption{The partially ordered set of minimal paths when $m=4$.}
\label{graph of paths}
\end{figure}

{\em Proof.}  Suppose that $K_{p/q}$ has 4 or less distinct boundary slopes and strongly positive continued fraction expansion $p/q = 0 + [a_1,a_2,..., a_m]$.  By Lemma~\ref{bounds on number of paths}, it follows that  $m\le 4$.    We will examine only the case $m=3$ and leave the remaining cases, which are handled in a similar manner,  to the interested reader.  

 Figure~\ref{graph of paths} shows all  possible minimal paths for $p/q=[a_1, a_2, a_3]$ arranged according to the partial ordering on paths.  Each of the five figures is made up of three main triangles which are subdivided into $a_1$, $a_2$, and $a_3$ smaller triangles, respectively, as we move from left to right.  The five paths $\gamma_\text{upper}$, $\gamma_1$, $\gamma_2$, $\gamma_3$, and $\gamma_\text{lower}$ are depicted with heavy dark lines.  Also listed is the unadjusted slope $m$ for each path under the path's label.

 Since $a_1 >1$ and $a_3 > 1$ we have that $\gamma_\text{upper}$, $\gamma_1$, $\gamma_2$, and $\gamma_\text{lower}$ are minimal. If $\gamma_3$ is minimal and  $a_1 \neq a_3$, then we would have 5 distinct boundary slopes.  Therefore, we have two cases to consider.

Assume $\gamma_3$ is not minimal, which requires $a_2 =1$, and also that $a_1\ne a_3$.   Thus, $p/q = 0 + [a_1,1,a_3]$ and we have four distinct boundary slopes.   In order to compute the (adjusted) slopes we need to determine the unique even path.  Notice that 
$$\frac{p}{q}=[a_1,1,a_3] = \frac{1+a_3}{a_1 + a_3 +a_1 a_3}$$
and so in order for $q$ to be odd we need at least one of $a_1$ or $a_3$ to be odd. This gives three sub-cases which are listed below.

\begin{center}
\begin{tabular}{|l|l|l|l|}
\hline
$a_1$ & $a_3$ & even path & boundary slopes \\
\hline
even & odd & $\gamma_1$ & $\{-2a_1,0,-2a_1+2a_3,2a_3+2\}$ \\
odd & even & $\gamma_2$ & $\{-2a_3,2a_1-2a_3,0,2a_1+2\}$ \\
odd & odd & $\gamma_\text{lower}$ & $\{-2a_1-2a_3-2,-2a_3-2,-2a_1-2,0\}$ \\ \hline
\end{tabular}
\end{center}

These slopes appear in parts $4_{ii}$ and $4_{iii}$ of Theorem~\ref{list of knots with given number of slopes}.  The first two rows correspond to equivalent knots and only appear in one row of the theorem. 

Next assume that $\gamma_3$ is not minimal and $a_1 = a_3$.  In this case, we have 3 distinct boundary slopes and
$$\frac{p}{q}=[a_1,1,a_1] = \frac{1+a_1}{a_1^2 + 2 a_1}.$$
Here $q$ is odd if and only if $a_1$ is odd. In this case, the even path is $\gamma_\text{lower}$ and the boundary slopes are
$$\{-4a_1-2,-2a_1-2,0\}$$
where the second slope $-2a_1-2$ has multiplicity two, that is, it corresponds to two paths,  $\gamma_1$ and $\gamma_2$.  This data appears in part $3_{iii}$ of Theorem~\ref{list of knots with given number of slopes}.

Finally, assume that $\gamma_3$ is minimal, which means $a_2>1$, and $a_1=a_3$. We now have four distinct boundary slopes and 
$$\frac{p}{q}=[a_1,a_2,a_1] = \frac{1+a_1 a_2}{a_1^2 a_2 + 2 a_1}.$$
In order for $q = a_1^2 a_2 + 2 a_1$ to be odd we must have that both $a_1$ and $a_2$ are odd.  This implies that $\gamma_\text{lower}$ is the even path and we obtain the following boundary slopes:
$$\{-4a_1-2a_2, -2a_1-2a_2, -2a_2, 0\}.$$
Notice that the second boundary slope $-2a_1-2a_2$ has multiplicity two.
This data appears in part $4_i$ of Theorem~\ref{list of knots with given number of slopes}.  The remaining cases follow by a similar analysis with $m=1$, $m=2$, or $m=4$. \hfill $\Box$

By working out the fractions $p/q$ for each case of Theorem~\ref{list of knots with given number of slopes} we obtain the following corollary which strengthens Theorem~3 of \cite{MMR:2008}. As in Theorem~\ref{list of knots with given number of slopes}, we consider a knot and its mirror image to be equivalent.

\begin{corollary} Let $K$ be a 2-bridge knot.

\begin{itemize}
\item[i.] $K$ has exactly two distinct boundary slopes if and only if $K \equiv K_{1/q}$.
\item[ii.] $K$ has exactly three distinct boundary slopes if and only if $K$ is not equivalent to a knot from part (i) and $K \equiv K_{p/q}$ for which either 
\begin{itemize}
\item[a.] $p  \mid q-1$ or 
\item [b.] $p^2 = q+1$.
\end{itemize}
\item[iii.] $K$ has exactly four distinct boundary slopes if and only if $K$ is not equivalent to a knot from parts (i) or (ii) and $K \equiv K_{p/q}$ for which either 
\begin{itemize}
\item[a.] $p + 1  \mid q$ and $q \mid p^2-1$, 
\item[b.] $p \mid q+1$,  
\item[c.] $(p-1)^3 = q^2$, or  
\item[d.] $p/q = 8/21$.
\end{itemize}
\end{itemize}
\end{corollary}

\section{Ohtsuki, Riley, Sakuma  knot pairs}

\begin{figure}[t]
\psfrag{X}{$K_{p/q}$}
\psfrag{Y}{$\beta$}
\psfrag{Z}{$K_{p'/q'}$}
\psfrag{a}{$a_1$}
\psfrag{b}{$a_2$}
\psfrag{c}{$a_3$}
\psfrag{g}{$2c_1$}
\psfrag{h}{$2c_2$}
\psfrag{d}{$\beta$}
\psfrag{e}{$\beta^{-1}$}
\psfrag{f}{$\beta_{-}$}
\psfrag{j}{$f$}
\psfrag{k}{$S^2 \times I$}
\psfrag{i}{$B^3$}

     \begin{center}
     \leavevmode
     \scalebox{.75}{\includegraphics{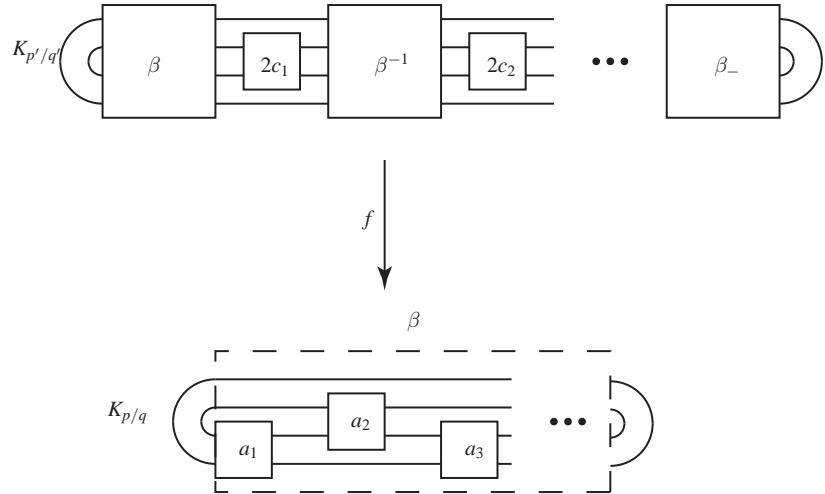}}
     \end{center}
\caption{The Ohtsuki, Riley, Sakuma construction.}
\label{ORS construction}
\end{figure}

In this section we briefly review the method of Ohtsuki, Riley, and Sakuma for constructing a pair of two bridge knots $K_{p'/q'}$ and $K_{p/q}$ with $K_{p'/q'} \ge K_{p/q}$. The reader is encouraged to consult \cite{ORS} for a more detailed description. One begins with any 2-bridge knot (or link) $K_{p/q}$ given by the four-plat shown in Figure~\ref{ORS construction}. The 4-string braid defining $K_{p/q}$ is denoted by $\beta$.  There are three associated braids $\beta^{-1}, \beta_{-}$, and $\beta^{-1}_{-}$ which are obtained by rotating $\beta$ 180 degrees around an axis perpendicular to the plane of the diagram, reflecting $\beta$ through a plane perpendicular to the plane of the diagram, and the composition of these two motions, respectively.  In order to construct a 2-bridge knot $K_{p'/q'}$ that is greater than $K_{p/q}$ we consider a 4-plat with an odd number of ``boxes'' each containing either $\beta, \beta^{-1}, \beta_{-}$, or $\beta^{-1}_{-}$.  Starting from the left in Figure~\ref{ORS construction}, the first box contains $\beta$ and every other box after that contains $\beta$ or $\beta_{-}$. The remaining boxes each contain $\beta^{-1}$ or $\beta^{-1}_{-}$. Also between each braid box in $K_{p'/q'}$ we may insert an even number of half-twists in the middle two strands of the 4-plat.  Finally, a {\it branched fold map} $f$ is constructed from the complement of $K_{p'/q'}$ onto the complement of $K_{p/q}$ as follows.  The complement of $K_{p'/q'}$ is cut along a collection of parallel 2-spheres into two 3-balls and a number of $S^2 \times I$'s. The complement of $K_{p/q}$ is decomposed into two 3-balls and a single $S^2 \times I$.  A continuous mapping is then defined in a piecewise manner.  First each 3-ball ``upstairs'' is mapped by the identity onto a corresponding 3-ball ``downstairs.'' Next, each $S^2 \times I$ containing one braid box upstairs is mapped homeomorphically onto the $S^2 \times I$ downstairs in a way that depends on the presence of $\beta, \beta^{-1}, \beta_{-}$, or $\beta^{-1}_{-}$.  Finally, the mapping is extended to the remaining $S^2 \times I$'s upstairs onto the 3-balls downstairs using 2-fold branched mappings that depend on the combinations of  $\beta, \beta^{-1}, \beta_{-}$, or $\beta^{-1}_{-}$ in the adjacent components.  From the way that $f$ is constructed, the two meridional generators associated with the bridges at either end of the 4-plat diagram upstairs are taken to the corresponding generators downstairs. Furthermore, the longitude upstairs is taken to a power of the longitude downstairs. Thus $f$ induces an epimorphism on fundamental groups which preserves peripheral structure.

Given a 2-bridge knot $K_{p/q}$, Ohtsuki, Riley and Sakuma determine which knots $K_{p'/q'}$ are produced by their construction. In particular, they show that $p'/q'$ gives such a knot if and only if it lies in the orbit of $p/q$ under the action of $\Gamma_\infty * \Gamma_{p/q}$, where $\Gamma_{r/s}$ is the infinite dihedral group generated by all reflections in edges of the Farey graph which end at $r/s$. Let $A_1$ and $A_2$ be generators of $\Gamma_\infty$ defined by  reflection in the edges $<1/0, 0/1>$ and $<1/0, 1/1>$, respectively. For $p/q=[a_1, a_2, \dots, a_m]$ let $B_1$ and $B_2$ be generators of $\Gamma_{p/q}$ defined by reflection in the edges $<p/q, [a_1, a_2, \dots, a_{m-1}]>$ and $B_2=<p/q, [a_1, a_2, \dots, a_m-1]>$. Thus $p'/q'=W(p/q)$ where $W$ is a word in $A_1, A_2, B_1$ and $B_2$. Since $B_1$ and $B_2$ fix $p/q$, and $A_1$ and $A_2$ take $p/q$ to a fraction representing the same knot, we may assume that we begin by applying a nontrivial word in $A_1$ and $A_2$ and end by applying a nontrivial word in $B_1$ and $B_2$. Hence we may write $W$ as 
\begin{equation}\label{W word}
W=W_1 W_2\dots W_n,  \mbox{ with $n$ even},
\end{equation}
where
\begin{equation}\label{W syllables} 
W_i=\left\{
\begin{array}{lll}
(B_1B_2)^{c_i}B_1^\frac{1-\eta_i}{2}, &c_i \in \mathbb Z,  \eta_i \in \{1, -1\}, &\mbox{ if $i$ is odd, and}\\ 
&&\\
(A_1A_2)^{c_i}A_1^\frac{1-\eta_i}{2}, &c_i \in \mathbb Z,  \eta_i \in \{1, -1\},& \mbox{ if $i$ is even.}
\end{array}
\right.
\end{equation}
 Ohtsuki, Riley and Sakuma then show (Lemma 5.3 in \cite{ORS}) that
\begin{equation}
\label{ORS formula}
p'/q'=[ \epsilon_1 {\bf a}, 2 \epsilon_1 c_1, \epsilon_2 {\bf a}^{-1}, 2 \epsilon_2 c_2, \epsilon_3 {\bf a}, 2 \epsilon_3 c_3, \epsilon_4 {\bf a}^{-1},\dots, \epsilon_{n+1} {\bf a}],
\end{equation}
where ${\bf a}$ is the vector ${\bf a}=(a_1, a_2, \dots, a_m),\  {\bf a}^{-1}$ means ${\bf a}$ written backwards, $\epsilon_1 = 1$, and $\epsilon_i=\prod_{j=1}^{i-1}-\eta_j$.
The following lemma will prove helpful in simplifying continued fractions. 

\begin{lemma}
\label{cf moves} Let $\bf a$ and $\bf b$ be vectors of integers, possibly empty, and $m$ and $n$ be any integers.  Then
\begin{enumerate}
\item[i.] $r+[{\bf a}, m, -n, {\bf b}] =r+ [{\bf a}, m-1,1,n-1, -{\bf b}]$ and
\item[ii.] $r+[{\bf a}, m, 0, n, {\bf b}] = r+[{\bf a}, m+n, {\bf b}]$
\end{enumerate}
\end{lemma}
 
\noindent {\em Proof.} 
If $p/q=r_0+[r_1, r_2, \dots, r_k]$ then it is well known that
$$\pm \begin{pmatrix}p\\q\end{pmatrix}= M_{(r_0, r_1, \dots, r_k)}\begin{pmatrix}1\\0\end{pmatrix},$$
where
$$\ M_{(r_0, r_1, \dots, r_k)}= 
\begin{pmatrix}r_0&1\\1&0\end{pmatrix}
\begin{pmatrix}r_1&1\\1&0\end{pmatrix}\cdots
\begin{pmatrix}r_k&1\\1&0\end{pmatrix}.$$
It is an easy exercise to show that
$$\begin{pmatrix}-1&0\\0&1\end{pmatrix} M_{\bf b}=(-I)^{|{\bf b}|}M_{-{\bf b}}\begin{pmatrix}-1&0\\0&1\end{pmatrix},$$
where $|{\bf b}|$ is the dimension of the vector ${\bf b}$ and $I$ is the identity matrix.
The proof of Lemma~\ref{cf moves}({\it i}) follows from

\begin{eqnarray*}
M_{{\bf a}, m, -n, {\bf b}} \begin{pmatrix}1\\0\end{pmatrix}&=&
M_{\bf a}\begin{pmatrix}m&1\\1&0\end{pmatrix} \begin{pmatrix}-n&1\\1&0\end{pmatrix} M_{\bf  b}
\begin{pmatrix}1\\0\end{pmatrix}\\
&=&
M_{\bf a}\begin{pmatrix}m-1&1\\1&0\end{pmatrix} 
\begin{pmatrix}1&1\\1&0\end{pmatrix}
\begin{pmatrix}n-1&1\\1&0\end{pmatrix}
\begin{pmatrix}-1&0\\0&1\end{pmatrix}
M_{\bf b}
\begin{pmatrix}1\\0\end{pmatrix}\\
&=&
(-I)^{|{\bf b}|} M_{{\bf a}, m-1, 1, n-1, -{\bf b}}
\begin{pmatrix}-1&0\\0&1\end{pmatrix}
\begin{pmatrix}1\\0\end{pmatrix}\\
&=&
(-I)^{|{\bf b}|} M_{{\bf a}, m-1, 1, n-1, -{\bf b}}
\begin{pmatrix}-1\\0\end{pmatrix}\\
&=&
(-I)^{|{\bf b}|+1} M_{{\bf a}, m-1, 1, n-1, -{\bf b}}
\begin{pmatrix}1\\0\end{pmatrix}
\end{eqnarray*}

The second part of the lemma follows in a similar way.
 \hfill $\Box$
 
We are now prepared to show that in any ORS pair $(K_{p'/q'}, K_{p/q})$, if $K_{p'/q'}$ is not a torus knot, then it has at least five distinct boundary slopes.

\begin{theorem}
\label{5slopes}
If $K_{p'/q'} \ge K_{p/q}$ is a nontrivial ORS pair, then either
\begin{enumerate}
 \item[{\it i}.] $K_{p'/q'}$ and $K_{p/q}$ are both torus knots and $K_{p'/q'}$ has precisely two distinct boundary slopes, or
 \item [{\it ii}.]$K_{p'/q'}$ has at least five distinct boundary slopes.
 \end{enumerate}
 \end{theorem}

 \noindent {\em Proof.} Let  $p/q=[a_1, a_2, \dots, a_m]=[{\bf a}]$ be a strongly positive continued fraction.  Suppose that $p'/q' = W(p/q)$ where $W$ is described in Equations (\ref{W word}) and (\ref{W syllables}).  By Equation~(\ref{ORS formula}) we have that
\begin{equation}
\label{ORS formula 2}
p'/q'=[ \epsilon_1 {\bf a}, 2 \epsilon_1 c_1, \epsilon_2 {\bf a}^{-1}, 2 \epsilon_2 c_2, \epsilon_3 {\bf a}, 2 \epsilon_3 c_3, \epsilon_4 {\bf a}^{-1},\dots, \epsilon_{n+1} {\bf a}].
\end{equation}
 
Suppose $c_i=0$ for some $i$. Because $W_i$ is not the identity, it follows that $\eta_i=-1$ and hence $\epsilon_i=\epsilon_{i+1}$. Therefore, on either side of $2 \epsilon_i c_i=0$, the continued fraction appears as
$p'/q'=[\dots,\epsilon_i {\bf a}^{\pm 1}, 0, \epsilon_i {\bf a}^{\mp 1},\dots]$. Using Lemma~\ref{cf moves}, we may eliminate the zero to combine $\epsilon_i {\bf a}^{\pm 1}$ and $\epsilon_i {\bf a}^{\mp 1}$ into a single strongly positive or negative vector. 

Let $1\le j_1<j_2<\dots <j_k\le n$ be the indices of the nonzero $c_i$'s.
Eliminating each of the $n-k$ zero entries as just described gives the following continued fraction expansion with $n(m-1)+m+2k$ partial quotients,
\begin{equation}
\label{nozerosfrac}
p'/q'=[ \epsilon_1 {\bf v}_1, 2 \epsilon_{j_1} c_{j_1}, \epsilon_{j_2} {\bf v}_2, 2 \epsilon_{j_2} c_{j_2},\epsilon_{j_3} {\bf v}_3,\dots, 2 \epsilon_{j_k} c_{j_k}, \epsilon_{n+1} {\bf v}_{k+1}]
\end{equation}
where $\epsilon_1=1$,  each ${\bf v}_i$ is a strongly positive vector,  and $c_{j_i} \neq 0$ for all $i$.  Notice that if $m=1$, then ${\bf a}^{-1}={\bf a}=a_1$ and each vector ${\bf v}_i$ is a single integer given by ${\bf v}_i=(j_i-j_{i-1})a_1$ if $i\le k$ and ${\bf v}_{k+1}=(n+1-j_k)a_1$.  Otherwise, if $m>1$,  $|{\bf v}_i| \ge 2$ for all $i$. 

We now claim that the continued fraction given in Equation~(\ref{nozerosfrac}) can be changed into a strongly positive one with at least as many partial quotients.  Before proving this claim, consider its consequences. If $m=1$ then the number of partial quotients is at least $2k+1$ which can only be less than 6 if $k=0, 1$, or $2$. If $m>1$, the number of partial quotients can only be less than 6 if $n=m=2$ and $k=0$. Hence, except in the cases where $m=1$ and $0\le k\le 2$, or $n=m=2$ and $k=0$, we can express $p'/q'$ as a strongly positive continued fractions with at least 6 partial quotients. Lemma~\ref{bounds on number of paths} now implies that $K_{p'/q'}$ has at least 5 distinct boundary slopes. Thus, after proving this claim, we must analyze these four cases.

To prove the claim, we first assume that $m>1$ and proceed by induction on the number of sign changes in the sequence of partial quotients in Equation~(\ref{nozerosfrac}). If there are none, then it is easy to check that Equation~(\ref{nozerosfrac}) is already strongly positive. If there are sign changes, then consider the smallest value of $r$ so that $\epsilon_{j_r}=1$ and either $ c_{j_r}$ or $\epsilon_{j_{r+1}}$ is negative. 

If $ c_{j_r}>0$ and  $\epsilon_{j_{r+1}}=-1$, then using Lemma~\ref{cf moves}, we obtain
\begin{eqnarray*}
p'/q'&=&[\dots, {\bf v}_r, 2  c_{j_r}, -{\bf v}_{r+1}, -2 c_{j_{r+1}},\dots]\\
&=&[\dots, {\bf v}_r, 2 c_{j_r}, -v_{r+1}^1, -v_{r+1}^2, \dots, -v_{r+1}^{|{\bf v}_{r+1}|}, -2 c_{j_{r+1}},\dots]\\
&=&[\dots,  {\bf v}_r, 2  c_{j_r}-1,1, v_{r+1}^1-1, v_{r+1}^2, \dots, v_{r+1}^{|{\bf v}_{r+1}|}, 2 c_{j_{r+1}},\dots].
\end{eqnarray*}
This continued fraction is now of the same form as the original, has fewer sign changes, and has more partial quotients. The result now follows by induction. It is important here to remember that since $m>1$, $|{\bf v}_i| \ge 2$ for all $i$.

If $ c_{j_r}<0$ and  $\epsilon_{j_{r+1}}=-1$, then using Lemma~\ref{cf moves}, we obtain
\begin{eqnarray*}
p'/q'&=&[\dots, {\bf v}_r, 2  c_{j_r}, -{\bf v}_{r+1}, -2 c_{j_{r+1}},\dots]\\
&=&[\dots, v_r^1, v_r^2, \dots, v_r^{|{\bf v}_r|}, 2  c_{j_r}, -{\bf v}_{r+1}, -2 c_{j_{r+1}},\dots]\\
&=&[\dots, v_r^1, v_r^2, \dots, v_r^{|{\bf v}_r|}-1, 1, -2  c_{j_r}-1, {\bf v}_{r+1}, 2  c_{j_{r+1}},\dots].
\end{eqnarray*}
Again, this continued fraction is now of the same form as the original, has fewer sign changes, and has more partial quotients. The result now follows by induction.

Finally, suppose that $ c_{j_r}<0$ and  $\epsilon_{j_{r+1}}=1$. Applying Lemma~\ref{cf moves} twice, we obtain
\begin{eqnarray*}
p'/q'&=&[\dots, {\bf v}_r, 2  c_{j_r}, {\bf v}_{r+1}, \dots]\\
&=&[\dots, v_1^1, \dots, v_r^{|{\bf v}_r|}-1, 1, -2c_{j_r}-2, 1, v_{r+1}^1-1, v_{r+1}^2, \dots, v_{r+1}^{|{\bf v}_{r+1}|}, \dots].
\end{eqnarray*}
If $c_{j_r}\ne -1$, then once again we have arrived a continued fraction of the same form as the original, but with fewer sign changes and more partial quotients. If instead, $c_{j_r}= -1$, then we may use Lemma~\ref{cf moves} to eliminate the zero entry and arrive at a continued fraction with the same form as the original, fewer sign changes, but now an equal number of partial quotients. In either case, the result follows by induction.

We now consider the case where $m=1$. Recall that if $m=1$, then ${\bf a}^{-1}={\bf a}=a_1$, ${\bf v}_i=(j_i-j_{i-1})a_1$ for $i\le k$, and ${\bf v}_{k+1}=(n+1-j_k)a_1$. Equation~(\ref{nozerosfrac}) is now of the form
\begin{equation}
\label{m=1 form}
p'/q'=[\epsilon_1 j_1 a_1, 2\epsilon_{j_1}c_{j_1}, \epsilon_{j_2}(j_2-j_1)a_1, 2\epsilon_{j_2}c_{j_2}, 
\dots, 2 \epsilon_{j_k}c_{j_k}, \epsilon_{n+1}(n+1-j_k)a_1].
\end{equation}
Because $a_1\ge 3$, this continued fraction is a member of a family $\cal F$ of continued fractions all having the form $[{\bf v}, {\bf w}]$ where 
\begin{itemize}
\item ${\bf v}$ is a strongly positive vector (possibly of length 1), 
\item the first entry of the vector ${\bf w}$ is negative, 
\item the entries of $\bf w$ alternate between being even and having magnitudes greater than 2 
(the first entry of $\bf w$ can have either property), and 
\item the last entry of $\bf w$ has magnitude greater than 2.
\end{itemize}
We call $\bf v$ the {\it strongly positive part} and ${\bf w}$ the {\it tail} of such a continued fraction. We shall prove, by induction on the length of the tail, that any such continued fraction can be changed to a strongly positive one with at least as many partial quotients. Since ${\bf v}$ is strongly positive, if ${\bf w}$ has length zero then we are done. Now consider such a continued fraction where ${\bf w}$ is not empty. We have
\begin{eqnarray*}
[\bf v, w]&=&[v_1, v_2, \dots, v_r, w_1, w_2, \dots, w_s]\\
&=&[v_1, v_2, \dots, v_r-1, 1, -w_1-1, -w_2, \dots, -w_s]
\end{eqnarray*}
If $-w_2>0$, or if $-w_2<0$ and $-w_1-1>1$, then we have arrived at a member of $\cal F$ with more partial quotients and a shorter tail. Hence the result follows by induction. Otherwise, we have $-w_1-1=1$ and $-w_2<0$. Thus
\begin{eqnarray*}
[\bf v, w]&=&[v_1, v_2, \dots, v_r, -2, w_2, \dots, w_s]\\
&=&[v_1, v_2, \dots, v_r-1, 1, 1, -w_2, \dots, -w_s]\\
&=&[v_1, v_2, \dots, v_r-1, 1, 0, 1, w_2-1, w_3, \dots, w_s]\\
&=&[v_1, v_2, \dots, v_r-1, 2, w_2-1, w_3, \dots, w_s].
\end{eqnarray*}
Because the entries of $\bf w$ alternate between being even and having magnitudes greater than 2, it must be the case that $w_2>2$ and hence $w_2-1>1$. Hence we have arrived at a member of $\cal F$ with the same number of partial quotients but with a shorter tail. Once again, the result now follows by induction.

We have now completed the proof of the claim, and it remains to analyze the four cases already mentioned above. If $m=n=2$ and $k=0$, we have $p/q=[a_1, a_2]$ and $p'/q'=[a_1, a_2, 0, a_2, a_1, 0, a_1, a_2]=[a_1, 2a_2, 2a_1, a_2]$. Comparing this strongly positive continued fraction to the table in Theorem~\ref{list of knots with given number of slopes} we see that $K_{p'/q'}$ cannot have less than 5 distinct boundary slopes.
 
If $m=1$ and $k=0$ the continued fraction for $p'/q'$ collapses to $p'/q'=[(n+1)a_1]$. Thus $K_{p'/q'}$ is a torus knot with exactly 2 slopes.
 
Next, suppose that $m=1$ and  $k=1$. Then exactly one of the $c_i$'s is nonzero, say $c_{j_1}$. Equation~\ref{nozerosfrac} now gives $p'/q'=[j_1 a_1, 2 \epsilon_{j_1}c_{j_1}, \epsilon_{n+1}(n+1-j_1)a_1]$. If all of these partial quotients are positive, then comparing this strongly positive continued fraction to the table in Theorem~\ref{list of knots with given number of slopes}, we see that $K_{p'/q'}$ cannot have less than  5 distinct boundary slopes.
If not, there are either one or two sign changes in the sequence of partial quotients. These considerations lead to the following cases which are handled in a fashion similar to that when $m>1$:
\begin{itemize}

\item  If $\epsilon_{j_1}c_{j_1}>0$ and $\epsilon_{n+1}=-1$, then
$p'/q'=[j_1a_1, 2\epsilon_{j_1}c_{j_1}-1, 1,  (n+1-j_1)a_1-1]$.
 
\item If $\epsilon_{j_1}c_{j_1}<0$ and $\epsilon_{n+1}=-1$, then
$p'/q'=[j_1a_1-1,1, -2\epsilon_{j_1}c_{j_1}-1, (n+1-j_1)a_1]$.

\item  If $\epsilon_{j_1}c_{j_1}=-1$ and $\epsilon_{n+1}=1$, then
$p'/q'=[j_1a_1-1, 2,  (n+1-j_1)a_1-1]$.

\item  If $\epsilon_{j_1}c_{j_1}<-1$ and $\epsilon_{n+1}=1$, then
$p'/q'=[j_1a_1-1,1, -2\epsilon_{j_1}c_{j_1}-2, 1,  (n+1-j_1)a_1-1]$.

\end{itemize}
In all four of these cases, comparison  to the table in Theorem~\ref{list of knots with given number of slopes} gives that $K_{p'/q'}$ cannot have less than  5 distinct boundary slopes.
 
Our final case to consider is when $m=1$ and $k=2$. Now two of the $c_i$'s are nonzero, say $c_{j_1}$ and $c_{j_2}$. Eliminating zeroes leads to
 $$p'/q'=[j_1 a_1, 2\epsilon_{j_1}c_{j_1}, \epsilon_{j_2}(j_2-j_1)a_1, 2\epsilon_{j_2}c_{j_2}, \epsilon_{n+1}(n+1-j_2)a_1].$$
If all of these partial quotients are positive, then comparing this strongly positive continued fraction to the table in Theorem~\ref{list of knots with given number of slopes}, we once again see that $K_{p'/q'}$ cannot have less than  5 distinct boundary slopes. If not, then using Lemma~\ref{cf moves} to move to a strongly positive expansion will increase the number of partial quotients to 6 or more, again resulting in 5 or more boundary slopes, unless the first entry of the tail is $-2$. In this case, that entry will remain even, and comparison with the table in Theorem~\ref{list of knots with given number of slopes} again shows that this knot cannot have less than 5 distinct boundary slopes.
\hfill $\Box$

The lower bound of five boundary slopes in part {\it ii} of Theorem~\ref{5slopes} is sharp.  For example, with $p/q=1/3=[3]$ and $p'/q'=7/45=[3,0,3,2,3]=[6,2,3]$ we have $K_{7/45} \ge K_{1/3}$. Furthermore, using Figure~\ref{graph of paths}, it is easy to show that $K_{7/45}$ has exactly 5 distinct boundary slopes.
 
Recall that if $K_{p'/q'}$ is a 2-bridge knot and $K_{p'/q'} \ge K_2$, then $K_2$ is also a 2-bridge knot.  Combining this fact with the result in Theorem~\ref{5slopes} gives the following corollary.
 
\begin{corollary}
If the answer to Question~\ref{ORS question} is yes, then 2-bridge knots with exactly 3 or 4 distinct boundary slopes are minimal with respect to the Silver-Whitten partial ordering.
\end{corollary}

Notice that if $K$ is a 2-bridge knot with exactly 2 distinct boundary slopes, then $K$ is the torus knot $T_{2, q}$ and is minimal if and only if $q$ is prime.

\section{Two-bridge knots with three distinct boundary slopes}

In this section we provide evidence for an affirmative answer to Question~\ref{ORS question} by proving that 2-bridge knots with exactly three distinct boundary slopes are in fact minimal.

\begin{theorem}
If $K$ is a 2-bridge knot with exactly 3 distinct boundary slopes, then $K$ is minimal with respect to the Silver--Whitten partial order.
\label{minknots}
\end{theorem}
{\em Proof.} Suppose $K_1$ is a 2-bridge knot with exactly three distinct boundary slopes and that $K_1\ge_d K_2$ with $K_2$ nontrivial. Then, by Theorem~\ref{properties of partial order} and Corollary~\ref{boundary slopes are subsets},  $K_2$ is a 2-bridge knot with exactly two or three distinct boundary slopes. If $K_2$ has two distinct slopes then it is a torus knot and so, by Gonz\'alez-Acu\~na and Ram\'irez  \cite{GR:2001}, $K_1$ is given by the ORS construction.  Now Theorem~\ref{5slopes} contradicts that $K_1$ has exactly three boundary slopes. Therefore, from Theorem~\ref{list of knots with given number of slopes} there are three possibilities for $K_1$ and three for $K_2$, giving a total of 9 different cases to consider. In what follows, the symbol $(3_{ii}, 3_{i})$, for example, will be used to denote the case where $K_1$ is the second type of knot with three slopes and $K_2$ is the third type of knot with three slopes, as listed in Theorem~\ref{list of knots with given number of slopes}.

Since boundary slopes of 2-bridge knots are always even, we first consider them mod 4. If $K_2$ has $r$ distinct slopes equal to 0 mod 4 and $s$ distinct slopes equal to 2 mod 4, then $K_1$ has the same number of slopes of each type because of Corollary~\ref{boundary slopes are subsets} and the fact that $d$ is odd. This observation rules out four cases: $(3_i, 3_{ii}), (3_i, 3_{iii}), (3_{ii}, 3_i), (3_{iii}, 3_i)$.

Before considering any other cases we note that the Alexander polynomial of of a 2-bridge knot with  three boundary slopes is given in the following table. Note that in each case, no nontrivial constant can be factored from each polynomial.

\begin{center}
\begin{tabular}{|ll|l|}
\hline
& knot type & Alexander polynomial\\
\hline
& & \\
$3_i$ & $0+[a_1, a_2],\ a_1$ even, $a_2$ even & $\displaystyle -\frac{a_1 a_2}{4}+ \left(1-\frac{a_1 a_2}{2} \right)t-\frac{a_1 a_2}{4}t^2$\\
& &  \\
$3_{ii}$ & $0+[a_1, a_2],\ a_1$ even, $a_2$ odd & $\displaystyle \frac{a_2+1}{2} + a_2 \sum_{i=1}^{a_1-1} (-1)^i  t^i + \frac{a_2+1}{2} t^{a_1}$\\
& & \\
$3_{iii}$ & $0+[a_1, -1, a_1],\ a_1$ odd&$\displaystyle \frac{(a_1+1)^2}{4}+ \left(1-\frac{(a_1+1)^2}{2}\right)t+\frac{(a_1+1)^2}{4}t^2$\\
& & \\
  \hline
 \end{tabular}
 \end{center}
 
\noindent{\bf Case $(3_i, 3_i)$}  Suppose $K_1$ corresponds to $0+[a_1,a_2]$ and $K_2$ corresponds to $0+[b_1, b_2]$. Because of Corollary~\ref{boundary slopes are subsets}, it follows that 
$$\{ -2db_1,0,2db_2\} = \{-2a_1,0,2a_2\}$$
and so $a_1 a_2=d^2 b_1 b_2$. Furthermore, $\Delta_{K_2}$ divides $\Delta_{K_1}$ which implies that  $d=\pm1$. Hence, $\{a_1,a_2\} = \{b_1,b_2\}$ and we have $K_2 \equiv K_1$.
 
\noindent {\bf Case $(3_{ii}, 3_{iii})$}  Suppose $K_1$ corresponds to $0+[a_1,a_2]$ and $K_2$ corresponds to $0+[b_1,1, b_1]$. From Corollary~\ref{boundary slopes are subsets} we obtain that $d<0$, $a_1=-(b_1+1)d$, and $a_1+a_2=-(2 b_1+1)d$. From these it follows that $a_2=-b_1 d$.  
Since the determinants of these knots are $a_1 a_2+1$ and $b_1^2+2 b_1$, respectively, it follows from Theorem~\ref{properties of partial order} that $k(b_1^2+2 b_1)=a_1a_2+1=b_1^2 d^2+b_1 d^2+1$ for some integer $k$. Hence $b_1$ divides 1 which is a contradiction.
 
\noindent {\bf Case $(3_{iii}, 3_{ii})$} Suppose $K_1$ corresponds to $0+[a_1,1,a_1]$ and $K_2$ corresponds to $0+[b_1,b_2]$. From Corollary~\ref{boundary slopes are subsets} we obtain that $d<0,\  a_1+1=-b_1 d$, and $2a_1+1=-(b_1+b_2)d$. From these we obtain that $d=-1,\ b_1=a_1+1$, and $b_2=a_1$. However, the division of the Alexander polynomials implies that $b_1=2$ and this contradicts that $a_1>1$.
 
 \noindent {\bf Case $(3_{iii}, 3_{iii})$} Suppose $K_1$ corresponds to $0+[a_1,1,a_1]$ and $K_2$ corresponds to $0+[b_1, 1, b_2]$. The division of the Alexander polynomials implies that $b_1=a_1$ and hence $K_2=K_1$.
  
\noindent {\bf Case $(3_{ii}, 3_{ii})$} Suppose $K_1$ corresponds to $0+[a_1,a_2]$ and $K_2$ corresponds to $0+[b_1, b_2]$. Using Corollary~\ref{boundary slopes are subsets} we immediately obtain that $a_1=b_1d$ and $a_2=b_2 d$. All evidence suggests that if $d \neq 1$, then the Alexander polynomials do not divide (in fact, they are almost certainly both irreducible), however, we were unable to prove this in general. Instead we turn to character varieties to settle this last case.
  
By Theorem~\ref{properties of partial order}, we have that $\phi^* : X(K_2) \rightarrow X(K_1)$ is an injective, algebraic  and closed (in the Zariski topology) mapping.  In particular, this implies that the image $\phi^*(X(K_2))$ is a subvariety of $X(K_1)$ that is birationally equivalent to $X(K_2)$.  By Theorem~6.5 of Macasieb, Petersen, and Van Luijk \cite{MPV:2009}, both $X(K_1)$ and $X(K_2)$ are irreducible curves and so we conclude that $X(K_1)$ and $X(K_2)$ are birationally equivalent, and therefore have the same genus.  Moreover, by Theorem~6.5 of \cite{MPV:2009}, the genus of $X(K_2)$ is
\begin{equation}
\label{genus2}
3 \left( \frac{b_2+1}{2} \right) \left( \frac{b_1}{2} \right) - \frac{b_2+1}{2} - 4 \frac{b_1}{2} +2.
\end{equation}
(To interpret Theorem~6.5 in our setting, note that $k=-b_2$ and $l = b_1$.)
Similarly, the genus of $X(K_1)$ is 
\begin{equation}
\label{genus1}
3 \left( \frac{d b_2+1}{2} \right) \left(\frac{d b_1}{2} \right)- \frac{d b_2+1}{2} - 4 \frac{d b_1}{2} +2.
\end{equation}
Equating (\ref{genus2}) and (\ref{genus1}) gives
$$d \left(3 d b_1 b_2 - 5b_1 - 2b_2 \right) = 3  b_1 b_2 - 5b_1 - 2b_2.$$
Since $d$, $b_1$, and $b_2$ are all positive, this can only be true if $d=1$, in which case $K_2 = K_1$.

\hfill $\Box$

Since Theorem~\ref{list of knots with given number of slopes} also classifies those 2-bridge knots with exactly four distinct boundary slopes, we could hope to apply arguments similar to those above to prove that 2-bridge knots with four slopes are minimal.  In particular, looking at boundary slopes will eliminate most of the 40 cases involved.  However, without additional information about the relevant Alexander polynomials or character varieties, we were unable to settle resolve several cases. We close with the following conjecture and question.

\begin{conjecture}
A 2-bridge knot with exactly four distinct boundary slopes is minimal with respect to the Silver-Whitten partial ordering.
\end{conjecture}

\begin{question}
Does there exist a non 2-bridge, non-minimal knot with exactly 3 distinct boundary slopes?
\end{question}


%
%
%

\bibliographystyle{gtart}
\bibliography{hs_11}


\end{document}